\documentclass[12pt]{article}

\usepackage{amsfonts,amsmath,amssymb,latexsym,theorem}

\newtheorem{corollary}{Corollary}
\newtheorem{theorem}[corollary]{Theorem}

\newtheorem{lemma}[corollary]{Lemma}

\theorembodyfont{\upshape}
\newtheorem{definition}[corollary]{Definition}
\newtheorem{remark}[corollary]{Remark}

\newtheorem{example}[corollary]{Example}

\newcommand{\pf}{{\em Proof\/:} }

\newcommand{\qed}{\nolinebreak\hfill$\Box$\vspace{3pt}}

\newcommand{\nz}{\mathbb N}
\newcommand{\gz}{\mathbb Z}

\begin{document}
\renewcommand{\textfraction}{0}
 
\title{Computing Symmetrized Weight Enumerators for Lifted Quadratic Residue Codes}
 \author{\normalsize
I.~M.~Duursma\\ 
\small Dept. of Mathematics\\[-5pt] 
\small University of Illinois\\[-5pt] 
\small Urbana-Champaign, IL\\[-5pt] 
\small duursma@math.uiuc.edu 
\and
\normalsize M.~Greferath\\
\small Dept. of Mathematics\\[-5pt] 
\small San Diego State University\\[-5pt]
\small San Diego, CA\\[-5pt] 
\small greferath@math.sdsu.edu}
\date{}
\maketitle
\thispagestyle{empty}
\begin{abstract}
The paper describes a method to determine symmetrized weight
enumerators of $\gz_{p^m}$-linear codes based on the notion of a
disjoint weight enumerator.
Symmetrized weight enumerators are given for the
lifted quadratic residue codes of length $24$ modulo $2^m$ 
and modulo $3^m$, for any positive $m$.
\end{abstract} 
\normalsize
 
\section{Introduction}
Ring-linear codes have gained importance since the beginning of the
last decade when it was discovered that certain non-linear binary
codes are actually linear over $\gz_4$, the ring of integers
modulo~$4$ (cf.~\cite{nechaev,hammons}). Since then a variety of
papers has appeared dealing with different foundational,
constructive and analytical aspects of ring-linear coding \cite{wood-dual,wood-weight,monomial}.
It appears that primary integer residue rings and more generally
Galois rings form a most important class of ring alphabets for
contemporary ring-linear coding theory \cite{babu-galois,byrne-galois,kuzmin-galois}. 
For ring-linear codes, the Hamming weight generalizes to the 
symmetrized weight which uses a partition of the alphabet into equivalence 
classes such that two letters are equivalent if and only if they are associated, meaning they are unit multiples of each other. In particular, for a $\gz_{p^m}-$linear code, 
the symmetrized weight uses a partition of the alphabet into $m+1$ 
equivalence classes and the symmetrized weight enumerators are 
polynomials in $m+1$ variables. 

The investigations in the current paper were motivated by the discovery of
some high-quality codes over $\gz_8$, $\gz_9$ and ${\rm
GR}(4,2)$.
The paper \cite{z8} introduces a binary
$(96,2^{37},24)$ code derived from a $\gz_8$-linear lift of the
extended binary Golay code. The properties of the binary code follow
from the symmetrized weight enumerator of the $\gz_8$-linear code.
Similarly, in \cite{z9} 
we describe a ternary $(72,3^{25},24)$-code
derived from a $\gz_9$-linear lift of the extended ternary QR-code 
of length $24$. Again, the properties of the ternary code are derived
from the symmetrized weight enumerator of the $\gz_9$-linear code.
Both the binary and the ternary code have more codewords than previously 
known codes of the same length and distance.

The investigation of these and other examples naturally
revealed the feasibility limitation of a brute-force determination
of symmetrized weight enumerators.
At the same time the quality
of the discovered codes suggested the development of a more rigid
and theoretical tool for the computation of structural parameters
for ring-linear codes.
This article therefore presents a method to efficiently compute
weight enumerators of linear codes over primary integer residue
rings. For the lifted QR-codes of length $24$ 
over $\gz_8$ and $\gz_9$, respectively, the method reproduces
weight enumerators that were obtained by brute-force computation 
in \cite{z8} and \cite{z9}, respectively. 

The main ingredient for our method is the notion of a 
{\em disjoint weight enumerator\/} $A_{s,t}$.
Given a family $(C_s)_{s\in \nz}$ of
$\gz_{p^s}$-linear codes with the property that $C_{s+1}$ is a lift
of $C_s$ for all $s\in \nz$, the weight enumerator $A_{s,t}$
contains combined information of $C_s$ and $C_t^\perp$. 
Furthermore we will use what we call {\em partial weight
enumerators\/} $D_{i,j}$ in order to ease the computation of the
$A_{s,t}$. It turns out that only a finite number of $D_{i,j}$
need to be computed in order to determine the entire collection
$(A_{s,t})_{s,t \in \nz}$.
This allows us to extend the results for codes over $\gz_{8}$ and 
$\gz_{9}$, and to describe the symmetrized weight enumerators for 
lifted QR-codes of length $24$ over $\gz_{2^m}$ and $\gz_{3^m}$,
for any $m$. For further properties of disjoint weight
enumerators and partial weight enumerators, we refer to \cite{dcc}.

\section{Disjoint weight enumerators}

Let $C_m$ be a code with coefficients in $\gz/p^m\gz$.
The composition of a codeword $c \in  C_m$ is  the $m+1$ tuple
${\rm wt}(c):=(n_0,n_1,\ldots,n_{m-1},n_\infty)$, where $n_i$ is the number of
coefficients in $c$ that have $p$-adic valuation $i$. 

\begin{definition} The symmetrized weight
enumerator of the code $C_m$ is the polynomial in $m+1$ variables
\[
A(x_0,x_1,\ldots,x_{m-1}; - ; z) \; = \;
\sum_{c \in C_m} x_0^{n_0} x_1^{n_1} \cdots x_{m-1}^{n_{m-1}} z^{n_\infty}
\]
\end{definition}

The variables in the weight enumerator are divided over three groups.
The second group of variables is reserved to give
information about the dual code. Thus,
\[
A(- ; y_0,y_1,\ldots,y_{m-1}; z) = \sum_{c \in C^\perp_m}
  y_0^{n_0} y_1^{n_1} \cdots y_{m-1}^{n_{m-1}} z^{n_\infty}
\]
will denote the dual weight enumerator.

\begin{example}
The Octacode (cf.~\cite{hammons}) is defined as an extended cyclic code of length $8$ over $\gz_4$. It is of rank $4$ over $\gz_4$ and has weight enumerator
\begin{eqnarray*}
A(u,v;-;z) &=& z^8 + 112z^3u^4v + 112zu^4v^3 + 16u^8 + 14z^4v^4 + v^8.
\end{eqnarray*}
The code is self-dual and the dual weight enumerator $A(-;u,v;z)$ equals $A(u,v;-;z)$. 
\end{example}

\begin{definition} \label{thm1}
For $m,\ell \geq 0$, we define {\em disjoint weight enumerators\/} inductively by
\begin{eqnarray} 
& &A_{m+1,\ell}(x_0,x_1,\ldots,x_{m}; y_0, y_1, \ldots, y_{\ell-1} ; z) \hspace{1cm}
\nonumber \\
& &\quad = \frac{1}{|C^\perp_1|} A_{m,\ell+1}(p x_0, p x_1,\ldots,p x_{m-1};
              y_0, y_1, \ldots, y_{\ell-1}, z-x_m ; z+(p-1)x_m).  \hspace{1cm} \label{eq:upA}
\end{eqnarray}
\end{definition}

The inverse transform 
\begin{eqnarray} 
& &A_{m,\ell+1}(x_0,x_1,\ldots,x_{m-1}; y_0, y_1, \ldots, y_{\ell} ; z) \hspace{1cm}
\nonumber \\
& &\quad = \frac{1}{|C|} A_{m+1,\ell}(x_0, x_1,\ldots, x_{m-1}, z-y_{\ell};
              p y_0, p y_1, \ldots, p y_{\ell-1} ; z+(p-1) y_{\ell}) \hspace{1cm} \label{eq:downA}
\end{eqnarray}
is easily verified. In particular the transform is symmetric in the two
groups of variables.

\begin{remark}
The above definition is compatible with the MacWilliams' transform 
for symmetrized weight enumerators over the alphabet $\gz_{p^m}$ that 
transforms $A_{m,0}$ into $A_{0,m}$. 
Thus the above transformations give a decomposition of the MacWilliams' transform 
into smaller steps, much like a Welsh-Hadamard transform. This is possible because
the $p$-adic weight function corresponds to a subgroup filtration of the alphabet 
$\gz_{p^m}$, such that at each step in the filtration the subgroups are
of same index $p$.
\end{remark}

\begin{example}
For a code over $Z_4$ the transforms in the defintion become
\begin{eqnarray*}
A_{2,0}(u,v;-;z) &=& A_{1,1}(2u;z-v;z+v) / |C^\perp_1| \\
A_{1,1}(u;v;z)   &=& A_{0,2}(-;v,z-u;z+u) / |C^\perp_1| 
\end{eqnarray*}
with inverse transforms
\begin{eqnarray*}
A_{1,1}(u;v;z)     &=& A_{2,0}(u,z-v;-;z+v) / |C_1| \\
A_{0,2}(-;u,v;z)   &=& A_{1,1}(z-v;2u;z+v) / |C_1| 
\end{eqnarray*}
In particular for the Octacode, 
\[
A_{1,1} = z^8 + 14u^4 z^4+ u^8 + 14v^4 z^4+ v^8 - 14u^4v^4
\]
\end{example}

\section{Partial weight enumerators}

The weight enumerator $A_{m+1,0}(x_0,x_1,\ldots,x_{m}; - ; z)$ of a
$\gz_{p^m}$-linear code $C$ contains a contribution 
$A_{m,0}(x_1,\ldots,x_{m}; - ; z)$ from the subcode $pC$. 
We introduce the partial weight enumerator $D_{m+1,0}$ that measures
the difference. So that
\begin{eqnarray} 
& &A_{m+1,0}(x_0,x_1,\ldots,x_{m}; - ; z) \hspace{1cm} \nonumber \\
& &\quad = A_{m,0}(x_1,\ldots,x_{m}; - ; z)
+ D_{m+1,0}(x_0,x_1,\ldots,x_{m}; - ; z).  \hspace{1cm}  \label{eq:A=A+D}
\end{eqnarray}
More generally we define an array of partial weight enumerators. 

\begin{definition} \label{thm2}
For $m,\ell\geq 0$ we define the {\em partial weight enumerators\/} $D_{m,\ell}$ as the sum of those monomials in $A_{m,\ell}$ that are divisible by $x_0 y_0$.
\end{definition}

Thus, the partial weight enumerators become the building blocks for the 
actual weight enumerators,
\begin{eqnarray}
& &A_{m,\ell}(x_0,\ldots,x_{m-1};y_0,\ldots,y_{\ell-1};z) \hspace{1cm} \nonumber \\
& &\quad = \sum_{i \leq m} \sum_{j \leq \ell}
D_{i,j}(x_{m-i},\ldots,x_{m-1};y_{\ell-j},\ldots,y_{\ell-1};z). \hspace{1cm} \label{eq:AsumD}
\end{eqnarray}


\begin{example} \label{ex:a11}
For the Octacode, the decomposition 
$A_{1,1} = D_{0,0}+D_{0,1}+D_{1,0}+D_{1,1}$ of $A_{1,1}(u;v;z)$ becomes
\[
\begin{array}{lcl}
D_{0,0}(-;-;z) = z^8, &  &D_{0,1}(-;v;z) = 14v^4 z^4+ v^8, \\
D_{1,0}(u;-;z) = 14u^4 z^4+ u^8, & &D_{1,1}(u;v;z) = -14u^4v^4.
\end{array}
\]
\end{example}

The terminology disjoint weight enumerator is motivated by the following theorem. 
It gives an
explicit expression for $D_{i,j}$ as a sum over pairs of codewords with disjoint support.

\begin{theorem}[see also \cite{dcc}] \label{thm:dcc}
Let $\{C_i: i \geq 0\}$ be a family of $\gz_{p^i}$-linear codes with the property that $C_{i+1}$ is a lift
of $C_i$ for all $i \geq 0$, and let $\{C_j^\perp: j \geq 0\}$ be the family of dual codes. For each $j \geq 0$,
there is a natural character family $\{\chi_v: v \in C_j^\perp\}$ such that
the partial weight enumerator $D_{i,j}$, for $i \geq 0$, is given by
\[
D_{i,j}(x_0,\ldots,x_{i-1};y_0,\ldots,y_{j-1};z)
  \; = \;  \sum_{\tiny\begin{array}{cc}
  u \in C_i,\, v \in C^\perp_j\\
  \bar{u}\neq 0\neq \bar{v}\\
  u \cap v = \emptyset
  \end{array}} \chi_v (u) {\bf x}^{{\rm wt}(u)} {\bf y}^{{\rm wt}(v)} z^{n-|u|-|v|}
\]
Here ${\bf x}^{{\rm wt}(u)} := x_0^{n_0} \cdots x_{i-1}^{n_{i-1}}$ where ${\rm wt}(u) = (n_0, \ldots, n_{i-1})$. 
An analogous definition is understood for ${\bf y}^{{\rm wt}(v)}$. And finally $z$ homogenizes the polynomial 
$D_{i,j}$ by counting the zeros outside $u \cup v$. 
\end{theorem}

The theorem shows that as $i,j$ increase we should expect fewer terms in the sum $D_{i,j}$
as the supports of codewords $u \in C_i$ and $v \in C^\perp_j$ are non-decreasing with $i$ and $j$.
In the next sections, it will turn out that this is the main reason that we will be able to 
compute and describe symmetrized weight enumerators efficiently. 

\begin{example}
The Octacode has empty sum $D_{2,1}$ since any word $u$ in $C_2$ with $\bar u \neq 0$ has Hamming weight
at least five and any word $v$ in $C^\perp_1$ with $\bar v \neq 0$ has Hamming weight at least four. And
the intersection of supports $v \cap u$ consists of at least one position. Thus we find
\begin{eqnarray*}
A_{2,1}(u,v;w;z) & = &D_{0,0}(-;-;z) + D_{1,0}(v;-;z) + D_{2,0}(u,v;-;z) + \\
                 & &  + \; D_{0,1}(-;w;z) + D_{1,1}(v;w;z) + D_{2,1}(u,v;w;z) \\
                 & = & (z^8) + (14v^4 z^4+ v^8) + (112z^3u^4v + 112zu^4v^3 + 16u^8)\\
                 & &  +\; (14w^4 z^4+ w^8) + (- 14v^4w^4) + (0).
\end{eqnarray*}
Consequently $A_{3,0}(u,v,w;-;z)$ can be computed via
\begin{eqnarray*}
A_{3,0}(u,v,w;-;z) & = & A_{2,1}(2u,2v;z-w;z+w)/2^{4} \\
& = & z^8 + 14w^4z^4 +  w^8+ 112z^3v^4w + 112zv^4w^3 + 16v^8\\
& & + \; 224z^3u^4v  + 672z^2u^4vw + 896zu^4v^3 + 672zu^4vw^2\\
& & + \; 256u^8 + 896u^4v^3w + 224u^4vw^3 + 16v^8
\end{eqnarray*}
In particular, the information to compute the $\gz_8$ weight enumerator is already 
contained in the $\gz_4$ weight enumerator. We say the Octacode has $p$-adic depth two.
\end{example}

\section{Computing symmetrized weight enumerators}

First we observe that the transform relations (\ref{eq:upA}) and (\ref{eq:downA}) for the weight
enumerators $A_{m,\ell}$ hold for the partial weight enumerators $D_{m,\ell}$ as well, 
provided both $m>0$ and $\ell>0$.

\begin{lemma} 
For $m,\ell > 0$,
\begin{eqnarray} 
& &D_{m+1,\ell}(x_0,x_1,\ldots,x_{m}; y_0, y_1, \ldots, y_{\ell-1} ; z) \hspace{1cm}
\nonumber \\
& &\quad = \frac{1}{|C^\perp_1|} D_{m,\ell+1}(p x_0, p x_1,\ldots,p x_{m-1};
              y_0, y_1, \ldots, y_{\ell-1}, z-x_m ; z+(p-1)x_m). \hspace{1cm} \label{eq:upD}
\end{eqnarray}
\begin{eqnarray} 
& &D_{m,\ell+1}(x_0,x_1,\ldots,x_{m-1}; y_0, y_1, \ldots, y_{\ell} ; z) \hspace{1cm}
\nonumber \\
& &\quad = \frac{1}{|C|} D_{m+1,\ell}(x_0, x_1,\ldots, x_{m-1}, z-y_{\ell};
              p y_0, p y_1, \ldots, p y_{\ell-1} ; z+(p-1) y_{\ell}).\hspace{1cm} \label{eq:downD}
\end{eqnarray}
\end{lemma}
\pf It suffices to consider (\ref{eq:upD}). With (\ref{eq:AsumD}) we obtain, for $\ell > 0$,
\begin{align*}
&A_{m+1,\ell}(x_0,x_1,\ldots,x_{m}; y_0, y_1, \ldots, y_{\ell-1} ; z)  \\
~=~&A_{m,\ell}(x_1,\ldots,x_{m}; y_0, y_1, \ldots, y_{\ell-1} ; z) +
A_{m+1,\ell-1}(x_0,x_1,\ldots,x_{m}; y_1, \ldots, y_{\ell-1} ; z)  \\
&-A_{m,\ell-1}(x_1,\ldots,x_{m}; y_1, \ldots, y_{\ell-1} ; z) +
D_{m+1,\ell}(x_0,x_1,\ldots,x_{m}; y_0, y_1, \ldots, y_{\ell-1} ; z).
\end{align*}
All terms but $D_{m+1,\ell}$ can be rewritten using the transform (\ref{eq:upA}), for $m > 0$.
But then also $D_{m+1,\ell}$ obeys the transform. \qed

For the cases $\ell=0$ or $m=0$, respectively, an extra term needs to be included 
in the right hand side of equation (\ref{eq:upD}) or (\ref{eq:downD}), respectively. 

\begin{lemma}
For $\ell=0$ and $m>0$, 
\begin{eqnarray}
& & D_{m+1,0}(x_0, \ldots, x_m ; - ; z) \nonumber \\
& &= \frac{1}{|C^\perp|} (D_{m,0}(p x_0, \ldots, p x_{m-1}; - ; z+(p-1)x_m) \hspace{1cm}
\nonumber \\
& &\quad \quad + D_{m,1}(p x_0, \ldots, p x_{m-1}; z-x_m ; z+(p-1)x_m) ) \hspace{1cm}
\label{eq:Dm+1}
\end{eqnarray}
For $m=0$ and $\ell>0$, 
\begin{eqnarray}
& & D_{0,\ell+1}(- ; y_0, \ldots, y_\ell ; - ; z) \nonumber \\
& &= \frac{1}{|C|} (D_{0,\ell}(- ; p y_0, \ldots, p y_{\ell-1}; - ; z+(p-1) y_\ell) \hspace{1cm}
\nonumber \\
& &\quad \quad + D_{1,\ell}(- ; p y_0, \ldots, p y_{\ell-1}; z-y_\ell ; z+(p-1) y_\ell) ) \hspace{1cm}
\label{eq:Dl+1}
\end{eqnarray}
\end{lemma}
\pf It suffices to consider (\ref{eq:Dm+1}). First we write 
\[
A_{m+1,0}(x_0, \cdots, x_m; - ; z)
= \frac{1}{|C^\perp_1|} A_{m,1}(p x_0, p x_1,\ldots,p x_{m-1};
              z-x_m ; z+(p-1)x_m).
\]
A decomposition similar to that used in the previous lemma applies to $A_{m,1}$, for $m > 0$, 
\begin{eqnarray*}
   & &A_{m,1}(p x_0, \ldots, p x_{m-1}; z-x_m ; z+(p-1)x_m) \\
 = & & A_{m-1,1}(p x_1, \ldots, p x_{m-1}; z-x_m ; z+(p-1)x_m) \\
   & &+A_{m,0}(p x_0, \ldots, p x_{m-1}; - ; z+(p-1)x_m) \\
   & &- A_{m-1,0}(p x_1, \ldots, p x_{m-1}; - ; z+(p-1)x_m) \\
   & &+ D_{m,1}(p x_0, \ldots, p x_{m-1}; z-x_m ; z+(p-1)x_m) 
\end{eqnarray*}
For the contribution $A_{m-1,1}$ we have
\[
A_{m,0}(x_1, \cdots, x_m; - ; z)
= \frac{1}{|C^\perp_1|} A_{m-1,1}(p x_1, \ldots, p x_{m-1};
              z-x_m ; z+(p-1)x_m).
\]
Collecting terms, using (\ref{eq:A=A+D}) twice, gives (\ref{eq:Dm+1}). 
The relation (\ref{eq:Dl+1}) follows by symmetry.
\qed

In general, we obtain the following recursive procedure to compute 
$A_{m+l+2,0}$ from $A_{m+l+1,0}$ for a known partial weight enumerator
$D_{m+1,l+1}$. 

\begin{theorem}[Computing symmetrized weight enumerators] \label{thm:comp}
Let the symmetrized weight enumerator $A_{m+\ell+1,0}$ be known together with a partial 
weight enumerator $D_{m+1,\ell+1}$, for some $m, \ell \geq 0$. Then $A_{m+\ell+2,0}$ is
obtained as follows
\begin{enumerate}
\item Use the transform (\ref{eq:upD}) repeatedly to obtain $D_{m+\ell+1,1}$ from $D_{m+1,\ell+1}$. 
\item Use the transform (\ref{eq:Dm+1}) to obtain $D_{m+\ell+2,0}$ from $D_{m+\ell+1,1}$ and $D_{m+\ell+1,0}$.
\item Use (\ref{eq:A=A+D}) to obtain $A_{m+\ell+2,0}$ from $A_{m+\ell+1,0}$ and $D_{m+\ell+2,0}$.
\end{enumerate}
\end{theorem}

\begin{remark}
To compute the symmetrized weight enumerators $A_{m,0}$ for all $m \geq 0$, it suffices 
with the theorem to give one partial weight enumerator $D_{i,j}$ for each value of $i+j \geq 0$.
The procedure succeeds in computing $A_{m,0}$ only in as far as we can determine $D_{i,j}$,
for $i+j \leq m$. As indicated after Theorem \ref{thm:dcc} this actually becomes easier as $i$ and $j$
increase as the $D_{i,j}$ tend to have fewer terms for larger values of $i$ and $j$. 
This behaviour of course is directly opposite to a brute-force computation, for which the 
complexity of computing $A_{m,0}$ grows exponentially with $m$. 
\end{remark}

\begin{remark}
The fact that the partial weight enumerators $D_{i,j}$ have few terms, especially for values
of $i$ and $j$ that are comparable, makes them the better choice for storing and representing 
the weight enumerator $A$. With the transforms used in the theorem they expand rapidly. 
\end{remark}

\begin{remark} \label{rem:comp}
The best way to compute the $D_{i,j}$ is recursively. Pairs of words $(u',v')$ that occur
in the sum $D_{i',j'}$, for large $i', j'$, reduce to pairs $(u,v)$ that occur in the sum 
$D_{i,j}$, for $i \leq i', j \leq j'$. Namely if $u', v'$ have disjoint support, then
certainly their reductions $u, v$ have disjoint support. Thus a critical stage is to
compute $D_{1,1}$, which runs over all pairs $(u,v) \in C \times C^\perp$ with disjoint
support. A large automorphism group will help to reduce the computations. For a given
complete list of disjoint words for the code and its dual over $\gz/p\gz$, it is relatively
straightforward to compute the $D_{i,j}$ as $i$ and $j$ increase. They can be computed
either directly, by lifting the words in the list, or by introducing unknowns for the
coefficients in $D_{i,j}$ and then computing just enough terms in the weight enumerator 
$A_{i+j,0}$ to determine the unknowns. The last option worked well for the quadratic
residue codes of length $24$, for which the results are given in the following sections.
\end{remark}

\section{Lifts of the $\bf[24,12,8]$ binary Golay code}

In this section we will apply theorem \ref{thm:comp} to 
the binary Golay code and verify results for the $\gz_4$ lifted code \cite{bonnsolecald95} and the 
$\gz_8$ lifted code \cite{z8}.
In fact we describe the symmetrized weight enumerator for any lift of the binary Golay code by
giving enough partial weight enumerators.
Computations as described in Remark \ref{rem:comp} were carried out in Magma \cite{Magma} and 
required less than a few minutes on a Pentium III 750MHZ laptop.
In comparison, computing the weight enumerator of just the $\gz/8\gz$ lifted Golay code via an exhaustive 
computation required more than a 24 hours, and is unfeasible for the larger alphabets.

\begin{theorem} \label{thm:qr24-2}
Let $\{ C_m : m \geq 0\}$ be the family of lifted extended cyclic quadratic residue codes
of length $24$ modulo $2^m$. The partial weight enumerators are
\begin{eqnarray*}
D_{0,0} &=&z^{24} \\[1ex]
D_{1,0} &=&759\,{z}^{16}{u}^{8}+2576\,{z}^{12}{u}^{12}
            +759\,{z}^{8}{u}^{16}+{u}^{24} \\[1ex]
D_{1,1} &=&759\,{u}^{16}{v}^{8}-2576\,{u}^{12}{v}^{12}+759\,{u}^{8}{v}^{16}
            -1518\,{z}^{8}{u}^{8}{v}^{8}  \\[1ex]
D_{2,1} &=&759 \cdot 16 \cdot u^{8}w^{8} \cdot 
         (-{v}^{2}{z}^{6}+2\,{v}^{4}{z}^{4}-{v}^{6}{z}^{2}) \\[1ex]
D_{2,2} &=&759 \cdot 64 \cdot u^{8}w^{8} \cdot 
         (-{v}^{2}{x}^{2}{z}^{4}+{v}^{4}{x}^{2}{z}^{2}
            +{v}^{2}{x}^{4}{z}^{2}+{v}^{4}{x}^{4}) \\[1ex]
D_{3,2} &=&759 \cdot 128 \cdot u^{8}x^{8} \cdot 
         (-{v}^{4}{y}^{2}{z}^{2}+2\,{v}^{4}w{y}^{2}z
            -{v}^{4}{w}^{2}{y}^{2}+{v}^{4}{y}^{4}-{v}^{2}{w}^{2}{y}^{4}) \\[1ex]
D_{3,3} &=&759 \cdot 128 \cdot u^{8}x^{8} \cdot 
         (-2\,{v}^{4}{y}^{4}+{v}^{4}{y}^{2}{t}^{2}
            +{v}^{2}{w}^{2}{y}^{4}-{v}^{2}{w}^{2}{y}^{2}{t}^{2}) \\[1ex]
D_{4,3} &=&759 \cdot 128 \cdot u^{8}y^{8} \cdot 
         (-{v}^{4}{t}^{2}{s}^{2}
            -{v}^{2}{w}^{2}{t}^{4}-2\,{v}^{2}wx{t}^{2}sz+2\,{v}^{2}w{x}^{2}{t}^{2}s) \\[1ex]
D_{4,4} &=&0, ~~\text{etc.}
\end{eqnarray*}
In each of the $D_{i,j}$'s, the variables $z, x_0, \ldots, x_{i-1}, y_0, \ldots, y_{j-1}$,
are replaced with the variables $z,u,v,w,x,y,t,s$, in that order.   
\end{theorem}

We look at two special cases.
The symmetrized weight enumerator $A_{2,0}(u,v;-;z)$ is particularly important. 
It was computed in \cite{bonnsolecald95}.
There it is shown that the theta series of the Leech lattice follows with the 
substitution
\[
z = \sum_{x \in 4 \gz} q^{x^2}, ~u = \sum_{x \in 4 \gz + 1} q^{x^2},
v = \sum_{x \in 4 \gz + 2} q^{x^2}
\]
Similar to Example \ref{ex:a11}, $A_{2,0}(u,v;-;z)$ has the following expression in 
terms of partial weight enumerators.

\begin{corollary}
For the quaternary extended Golay code,
\[
A_{2,0}(u,v;-;z) = A_{1,1}(2u;z-v;z+v) / 2^{12}
\]
\[
A_{1,1}(u;v;z) = D_{0,0}(-;-;z) + D_{1,0}(u;-;z) + D_{0,1}(-;v;z) + D_{1,1}(u;v;z) 
\]
\end{corollary}

Thus, for a known weight enumerator $A_{1,0}$ of the binary golay code, 
the only new information is given by the
coefficients $759, -2576, 759, -1518$ in $D_{1,1}$. They are uniquely determined with the
arguments from \cite{bonnsolecald95}. The first three coefficients are 
unique such that the Euclidean weights of the lifted code are divisible by eight. And 
the final coefficient is unique such that the lifted code contains no word with eight
units and sixteen zeros. In \cite{bonnsolecald95}, $A_{2,0}$ is first expressed in
terms of invariants before the coefficients can be determined. That part of the
computation is avoided in our computation. \\

Another important weight enumerator is $A_{3,0}(u,v,w;-;z)$ which yields the Hamming
weight enumerator for the non-linear binary code $(96,2^{36},24)$ used in \cite{z8} after
the substitution $(u,v,w,z) \mapsto (t^2, t^2, t^4, 1)$. 
We give the expression for $A_{3,0}(u,v;-;z)$ in terms of the partial weight enumerators.
The expanded expression covers one page \cite{z8}.

\begin{corollary}
For the $\gz_8$-linear extended Golay code,
\[
A_{3,0}(u,v,w;-;z)  =  A_{2,1}(2u,2v;z-w;z+w) / 2^{12}
\]
\begin{eqnarray*}
A_{2,1}(u,v;w;z) & = &A_{2,0}(u,v;-;z) + A_{1,1}(v;w;z) - A_{1,0}(v;-;z) + D_{2,1}(u,v;w;z) 
\end{eqnarray*}
\end{corollary}

\section{Lifts of the $\bf[24,12,9]$ ternary QR code}

We describe the symmetrized weight enumerator for any lift of the extended 
ternary quadratic residue code of length $24$.

\begin{theorem} \label{thm:qr24-3}
Let $\{ C_m : m \geq 0\}$ be the family of lifted extended cyclic quadratic residue codes
of length $24$ modulo $3^m$. The partial weight enumerators are
\begin{eqnarray*}
D_{0,0} &=&z^{24} \\[1ex]
D_{1,0} &=&4048\,{z}^{15}{u}^{9}+61824\,{z}^{12}{u}^{12}
             +242880\,{z}^{9}{u}^{15}+198352\,{z}^{6}{u}^{18}
             +24288\,{z}^{3}{u}^{21}+48\,{u}^{24} \\[1ex]
D_{1,1} &=&-16192\,{z}^{6}{u}^{9}{v}^{9}+12144\,{z}^{3}{u}^{9}{v}^{12}
+12144\,{z}^{3}{u}^{12}{v}^{9}-1104\,{u}^{12}{v}^{12} \\[1ex]
D_{2,1} &=&4048 \cdot 18 \cdot {u}^{9}{w}^{9} \cdot \left(
-{v}^{2}{z}^{4}+{v}^{3}{z}^{3}+{v}^{5}z-{v}^{6}-{u}^{3}{z}^{3}
+3\,{u}^{3}v{z}^{2}-3\,{u}^{3}{v}^{2}z+{u}^{3}{v}^{3}-{v}^{3}{w}^{3}
\right)  \\[1ex]
D_{2,2} &=&4048 \cdot 117 \cdot {u}^{9}w^{9} \cdot (-{v}^{3}{x}^{3}) \\[1ex]
D_{3,2} &=&4048 \cdot 27 \cdot {u}^{9}{x}^{9} \cdot \left(
-6\,{v}^{2}{y}^{2}wz+6\,{v}^{2}{y}^{2}{w}^{2}-{v}^{3}{y}^{3} \right) \\[1ex]
D_{3,3} &=&0, ~~\text{etc.}
\end{eqnarray*}
In each of the $D_{i,j}$'s, the variables $z, x_0, \ldots, x_{i-1}, y_0, \ldots, y_{j-1}$,
are replaced with the variables $z,u,v,w,x,y,t,s$, in that order.   
\end{theorem}

The weight enumerator $A_{2,0}(u,v;-;z)$ yields the Hamming
weight enumerator for the non-linear ternary code $(72,3^24,24)$ used in \cite{z9} after
the substitution $(u,v,z) \mapsto (t^2, t^3, 1)$. 
The expanded expression for $A_{2,0}(u,v;-;z)$ covers half a page \cite{z9}.

\begin{corollary}
For the $\gz_9$-linear extended QR-code of length $24$,
\[
A_{2,0}(u,v;-;z)  =  A_{1,1}(3u;z-w;z+2w) / 3^{12}
\]
\[
A_{1,1}(u;v;z) = D_{0,0}(-;-;z) + D_{1,0}(u;-;z) + D_{0,1}(-;v;z) + D_{1,1}(u;v;z) 
\]
\end{corollary}

\end{document}